\documentclass{amsart}

\newtheorem{thm}{Theorem}[section]

\newtheorem{lem}[thm]{Lemma}
\newtheorem{cor}[thm]{Corollary}

\theoremstyle{definition}

\theoremstyle{remark}

\newtheorem{remark}[thm]{Remark}

\numberwithin{equation}{section}

\begin{document}


\title{Best approximations for the weighted combination of the Cauchy--Szeg\"o kernel and its derivative in the mean}


\author{Viktor V. Savchuk}
\address{Institute of Mathematics of NAS of Ukraine, 
Kyiv, Ukraine, 01024}
\email{vicsavchuk@gmail.com}
\urladdr{https://www.researchgate.net/profile/Viktor-Savchuk} 

\author{Maryna V. Savchuk}
\address{National Technical University of Ukraine "Igor Sikorsky Kyiv Polytechnic Institute"}
\email{ma.savchuk@kpi.ua}
\urladdr{https://www.researchgate.net/profile/Mv-Savchuk}





\begin{abstract}
In this paper, we study an extremal problem involving best approximation in the Hardy space $H^1$ on the unit disk 
$\mathbb D$.  Specifically, we consider weighted combinations of the Cauchy–Szeg\"o kernel and its derivative, parameterized by an inner funtion $\varphi$ and a complex number  $\lambda$, and provide explicit formulas for the best approximation $e_{\varphi,z}(\lambda)$ by the subspace $H^1_0$. We also describe the extremal functions associated with this approximation. Our main result gives the form of $e_{\varphi,z}(\lambda)$ as a function of $\lambda$ and shows that, for the entire range of $\lambda$, the extremal function is linear in $\lambda$ and unique. We apply this result to establish a sharp inequality for holomorphic functions in the unit disk, leading to a new version of the Schwarz-Pick inequality.
\end{abstract}


\maketitle




\section{Introduction}

Let $\mathbb D:=\{z\in\mathbb C : |z|<1\}$, $\mathbb T:=\{t\in\mathbb T : |t|=1\}$ and $\mathrm{dm}$ is the normalised Lebesgue measure on $\mathbb T$.

The Hardy space $H^p$ for $1\le p\le\infty$ is the class of holomorphic in the unit disk $\mathbb D$ functions $f$ satisfied
\[
+\infty>\|f\|_p:=
\begin{cases}
	\displaystyle\sup_{0\le\rho<1}\int_\mathbb T|f(\rho t)|^p\mathrm{dm}(t),\hfill&\mbox{if}~ p>0,\cr\cr
	\displaystyle\sup_{z\in\mathbb D}|f(z)|,\hfill&\mbox{if}~p=1.
\end{cases}
\]

It is well known, that for each function $f\in H^p$, the nontangential limit, say $f(t),$ $t\in\mathbb T,$ exist almost everywhere on $\mathbb T$ and $f\in L^p(\mathbb T)$ with 
\[
\|f\|_p=\left(\int_\mathbb T|f|^p\mathrm{dm}\right)^{1/p}.
\]

Further, we denote by \(\|\cdot\|_p\) (\(1 \le p < \infty\)) the norm in both \(H^p\) and \(L^p(\mathbb{T})\).

The Cauchy--Szeg\"o kernel is the function
\[
C(t,z):=\frac{1}{1-t\overline z}
\]
defined on $\mathbb T\times\mathbb D$.

It is well known that Cauchy--Szeg\"o kernel have a reproducing property for the Hardy space $H^1$, that is
\[
f(z)=\int_\mathbb Tf(t)\overline{C(t,z)}\mathrm{dm}(t),\quad z\in\mathbb D,
\]
for every function $f\in H^1$.

For some function $\varphi\in L^1(\mathbb T)$, $\|\varphi\|_1>0$, complex numbers $\lambda\in\mathbb C$ and $z\in\mathbb D$,  by weighted combination of the Cauchy--Szeg\"o kernel and its derivative  we mean the function
\[
C_{\varphi,z,\lambda}(t):=\varphi(t)\frac{\partial }{\partial\overline z}C(t,z)+\lambda C(t,z),
\]
defined on $\mathbb T$.

Denote by
\[
e_{\varphi,z}(\lambda):=\inf\left\{\left\|\overline{C_{\varphi,z,\lambda}}+h\right\|_1 : h\in H^1_0\right\},
\]
where $H^1_0:=\{f\in H^1 : h(0)=0\}$, the best approximation of $\overline{C_{\varphi,z,\lambda}}$ in the mean on the unit circle $\mathbb T$ by Hardy subspace $H^1_0$.

The function $h\in H^1_0$ (if it exists) for which the last  infimum is attained is called the extremal function for $e_{\varphi,z}(\lambda)$. 

The study of best approximations of the Cauchy–Szeg\"o kernel, or its analogue for the real axis, by various functional subspaces is a fundamental extremal problem in both approximation theory and the theory of holomorphic functions \cite{Rogo}. Let us mention the articles \cite{Alp1,Alp2,Sav1,Sav2,Sav3,Sav4,Sav5,Sav6,Chaich,Sav7,Sav8,Sav9}, where this problem was solved in specific cases.

In particulary, in  \cite{Sav8}, it was shown that for $\varphi(t)=t^n$, $n\in\mathbb Z_+$,
\[
e_{\varphi,z}(0)=\frac{1}{1-|z|^2},~z\in\mathbb D,
\]
and the extremal function is unique, with $h\equiv 0$. For other choises of $\varphi$ and  $\lambda$, the value $e_{\varphi,z}(\lambda)$ remains unknown.

In this paper, we  give an explisit form of $e_{\varphi,z}(\lambda)$ as a function of $\lambda\in\mathbb C$ and discribe the corresponding extremal functions in the case when $\varphi$ is an inner function. 

Let us recall that an inner function is a bounded holomorphic function $\varphi$ in $\mathbb D$ whose radial limits
\[
\varphi^*(t)=\lim_{r\to 1-}\varphi(rt)
\]
satisfy the equality $|\varphi^*(t)|=1$ for almost all $t\in\mathbb T$.

\section{Main result}
The main result is following

\begin{thm}\label{Thm 1}
	Let $z\in\mathbb D$ and $\varphi$ is an inner function. Then for all $\lambda\in\mathbb C$ we have
	\begin{align}\label{eq1}
		e_{\varphi,z}(\lambda)
		&=
		\begin{cases}
			\displaystyle\frac{1}{1-|z|^2}+\frac{|\lambda|^2(1-|z|^2)}{4},&\hfill\mbox{if}\quad\displaystyle 0\le|\lambda|\le\frac{2}{1-|z|^2},\\
			|\lambda|,&\hfill\mbox{if}\quad\displaystyle |\lambda|\ge\frac{2}{1-|z|^2},
		\end{cases}
	\end{align}
	The extremal for $e_{\varphi,z}(\lambda)$ is only the function
	\begin{equation}\label{extremal g}
		h(t)=\beta^2\varphi(t)\frac{\partial}{\partial\overline z}C(t,z)+\lambda t\overline zC(t,z),
	\end{equation}
	where
	\[
	\beta=
	\begin{cases}
		\displaystyle\frac{\lambda(1-|z|^2)}{2},&\hfill\mbox{if}\quad \displaystyle|\lambda|\le\frac{2}{1-|z|^2},\\
		\mathrm e^{\mathrm i\arg\lambda},&\hfill\mbox{if}\quad \displaystyle|\lambda|\ge\frac{2}{1-|z|^2}.
	\end{cases}
	\]
\end{thm}

\begin{remark}
	It follows from \eqref{eq1} that 
	\[
	\lim_{|\lambda|\to+\infty}\frac{e_{\varphi,z}(\lambda)}{|\lambda|}=\min\left\{\left\|\overline{C(\cdot,z)}+h(\cdot)\right\|_1 : h\in H^1_0\right\}=1,
	\]
	where minimum is attained only for the function
	\[
	h(t)=t\overline zC(t,z).
	\]
\end{remark}

The proof of theorem 1 is dased on the following
\begin{lem}\label{Lem 1}
	Let $k\in L^1(\mathbb T)$. The function $h\in H^1_0$ is the unique best approximation to $k$ in $L^1(\mathbb T)$ i.e.
	\[
	\min\left\{\|k+h\|_1 : h\in H^1_0\right\}=\|k+h\|_1,	
	\]
	if and only if there exist an inner function $I$ and a non-negative function $p\in L^1(\mathbb T)$ such that  $k+h=\overline Ip$ for almost  all $t\in\mathbb T$. 
\end{lem}

A sufficient version of this assertion was formulated in \cite{Sav8}.

\begin{proof}[Proof of lemma~\ref{Lem 1}] Clearly, it is  enough to prove the lemma in the case $h\equiv 0$.
	
	$"\Rightarrow"$ By the duality theorem (see \cite{Garnett}, p. 129) there exist the unique inner function $I$ such that
	\[
	Ik=|k|
	\]
	a.e. on $\mathbb T$. This gives us the representation $k=\overline I|k|$ a.e. on $\mathbb T$.
	
	$"\Leftarrow"$ We have 
	\[
	I(k+g)=p+Ig
	\]
	for any $g\in H^1_0$ a.e. on $\mathbb T$. Since  $|k+g|=|p+Ig|$ and $Ig\in H^1_0$ a.e. on $\mathbb T$, we get
	\begin{align*}
		\inf\left\{\|k+g\|_1 : g\in H^1_0\right\}&=\inf\left\{\|p+Ig\|_1 : g\in H^1_0\right\}\\
		&=\inf\left\{\|p+g\|_1 : g\in H^1_0\right\}.
	\end{align*}
	However, $p\ge 0$ a.e. on $\mathbb T$. Thus by Theorem 1 in \cite{Sav4} the last infimum is attained only for the $g\equiv0$. So,
	\[
	\inf\left\{\|k+g\|_1 : g\in H^1_0\right\}=\|p\|_1=\|k\|_1
	\] 
	and the proof is complete.
\end{proof}

\begin{proof}[Proof of theorem~\ref{Thm 1}]	
	Take $\lambda$ such that $0\le |\lambda|\le 2/(1-|z|^2)$.
	Then the sum $\overline{C_{\varphi,z,\lambda}(t)}+h(t)$, where $t\in\mathbb T$ and the function $h$ is defined by  \eqref{extremal g},  with a little algebra take 
	the form 
	\begin{align*}
		\overline{C_{\varphi,z,\lambda}(t)}+h(t)&=\overline {t\varphi(t)}\left(\frac{1}{1-\overline tz}+\beta\frac{t\varphi(t)}{1-t\overline z}\right)^2\\
		&=\overline {t\varphi(t)}\left(\frac{1+\beta\varphi(t)\frac{t-z}{1-t\overline z}}{1-\overline tz}\right)^2\\
		&=\overline{\left(\frac{\varphi(t)\frac{t-z}{1-t\overline z}+\beta}{1+\beta \varphi(t)\frac{t-z}{1-t\overline z}}\right)}\left|\frac{1+\beta\varphi(t)\frac{t-z}{1-t\overline z}}{1-\overline tz}\right|^2.
	\end{align*}
	
	Since $0\le|\beta|\le 1$,  the first factor in the last equality is a complex conjugate of the inner function 
	\begin{equation}\label{inner funct}
		I(t)=\frac{\varphi(t)\frac{t-z}{1-t\overline z}+\beta}{1+\beta  \varphi(t)\frac{t-z}{1-t\overline z}}.
	\end{equation}
	
	Therefore, by lemma \ref{Lem 1}, $h$ is the unique function in $H^1_0$ that best approximates $\overline{C_{\varphi,z,\lambda}}$. Besides, by the Parseval identity, we get
	\begin{align*}
		\left\|\overline{C_{\varphi,z,\lambda}}+h\right\|_1&=\int_\mathbb T\left|\overline{C(t,z)}+\beta t\varphi(t)C(t,z)\right|^2\mathrm{dm}(t)\\
		&=\|C(\cdot,z)\|_2^2+|\beta|^2\|C(\cdot,z)\|_2^2\\
		&=\frac{1}{1-|z|^2}+|\beta|^2\frac{1}{1-|z|^2}.
	\end{align*}
	
	Now, let us assume that $|\lambda|\ge 2/(1-|z|^2)$. Then we have
	\begin{align*}
		\mathrm e^{-\mathrm i\arg\lambda}\left(\overline{C_{\varphi,z,\lambda}(t)}+h(t)\right)&=|\lambda|\frac{1-|z|^2}{|1-\overline tz|^2}+2\mathrm{Re}\frac{\mathrm e^{\mathrm i\arg\lambda}t\varphi(t)}{(1-t\overline z)^2}\\
		&=
		|\lambda|\frac{1-|z|^2}{|1-\overline tz|^2}+\frac{2}{|1-\overline tz|^2}\mathrm{Re}\left(\mathrm e^{\mathrm i\arg\lambda}\varphi(t)\frac{t-z}{1-t\overline z}\right)\\
		&\ge
		\frac{1}{|1-\overline tz|^2}\left(|\lambda|(1-|z|^2)-2\right)\\
		&\ge 0
	\end{align*}
	for all $t\in\mathbb T$. Therefore, by lemma \ref{Lem 1} $h$ defined by \eqref{extremal g} still  extremal also in this case. Besides by the Poisson formula, we get
	\begin{align*}
		\left\|\overline{C_{\varphi,z,\lambda}}+h\right\|_1&=\int_\mathbb T\left(|\lambda|\frac{1-|z|^2}{|1-\overline tz|^2}+2\mathrm{Re}\frac{\mathrm e^{\mathrm i\arg\lambda}t\varphi(t)}{(1-t\overline z)^2}\right)\mathrm{dm}(t)\\
		&=|\lambda|\int_\mathbb T\frac{1-|z|^2}{|1-\overline tz|^2}\mathrm{dm}(t)\\
		&=|\lambda|.
	\end{align*}
\end{proof}

We now demonstrate an application of Theorem \ref{Thm 1} to establish a sharp inequality for holomorphic functions in the unit disk, resulting in a new version of the Schwarz-Pick inequality.

\begin{cor}
	Let \( f \) be a holomorphic function in \( \mathbb{D} \) with \( \sup_{z \in \mathbb{D}} |f(z)| \leq 1 \), and let \( \varphi \) be an inner function such that \( f/\varphi \) is holomorphic in \( \mathbb{D} \). Then, for all \( z \in \mathbb{D} \) and \( \lambda \geq 0 \), we have
	\[
	\left| \left( \frac{f}{\varphi} \right)'(z) \right| + \lambda |f(z)| \leq \begin{cases}
		\displaystyle \frac{1}{1 - |z|^2} + \frac{\lambda^2 (1 - |z|^2)}{4}, & \hfill \text{if} \quad \displaystyle 0 \leq \lambda \leq \frac{2}{1 - |z|^2}, \\
		\lambda, & \hfill \text{if} \quad \displaystyle \lambda \geq \frac{2}{1 - |z|^2}.
	\end{cases}
	\]
	In particular,
	\begin{equation}\label{Schwarz-Pick}
		\left| \left( \frac{f}{\varphi} \right)'(z) \right| \leq \frac{1 - |f(z)|^2}{1 - |z|^2},~z\in\mathbb D.
	\end{equation}
	Equality is attained for the function given by \eqref{inner funct} with $\beta=0$.
\end{cor}

\begin{proof}[Proof of the corollary.]
	Assume for the moment that \( \lambda \) is complex. Then, by the Cauchy formula, we get
	\begin{align*}
		\left( \frac{f}{\varphi} \right)'(z) + \lambda f(z)&= \int_{\mathbb{T}} f(t) \overline{C_{\varphi,z,\lambda}(t)} \, \mathrm{dm}(t)\\
		&= \int_{\mathbb{T}} f(t) \left( \overline{C_{\varphi,z,\lambda}(t)} + h(t) \right) \mathrm{dm}(t), \quad z \in \mathbb{D},
	\end{align*}
	where \( h \in H^1_0 \) is arbitrary.
	
	By substituting \( \lambda = e^{i \theta} \lambda \), where \( \theta = \arg\left( \frac{f}{\varphi} \right)'(z) - \arg f(z) \) and \( \lambda \geq 0 \), it follows that
	\[
	\left| \left( \frac{f}{\varphi} \right)'(z) \right| + \lambda |f(z)| \leq \left\| \overline{C_{\varphi,z,\lambda}} + h \right\|_1.
	\]
	Taking the infimum over \( h \in H^1_0 \), we obtain the desired estimate by Theorem \ref{Thm 1}.
	
	To verify \eqref{Schwarz-Pick}, choose \( \lambda = \frac{2|f(z)|}{1 - |z|^2} \).
\end{proof}

\bibliographystyle{amsplain}

\end{document}